\documentclass[runningheads]{llncs}
\usepackage{graphicx}
\usepackage{tabularx,wrapfig,amsmath}
\begin{document}
\title{A short proof of the non-biplanarity of $K_9$\thanks{Supported by NSERC.}}
%
%
\author{Ahmad Biniaz}
\authorrunning{A. Biniaz}
%
\institute{School of Computer Science, University of Windsor\\
\email{ahmad.biniaz@gmail.com}}
\maketitle              
\begin{abstract}
Battle, Harary, and Kodama (1962) and independently Tutte (1963) proved that the complete graph with nine vertices is not biplanar. Aiming towards simplicity and brevity, in this note we provide a short proof of this claim.

\end{abstract}
\section{Introduction}

An embedding (or drawing) of a graph in the Euclidean plane is a mapping of its vertices to distinct points in the plane and its edges to smooth curves between their corresponding vertices. A planar embedding of a graph is a drawing of the graph such that no two edges cross. A graph that admits such a drawing is called planar. A {\em biplanar embedding} of a graph $H=(V,E)$ is a decomposition of $H$ into two planar graphs $H_1 = (V, E_1)$ and $H_2 = (V, E_2)$ such that $E_1 \cup E_2=E$ and $E_1 \cap E_2 = \emptyset$, together with planar embeddings of $H_1$ and $H_2$. In this case, $H$ is called {\em biplanar}. In other words, a graph is called biplanar if it is the union of two planar graphs; that is, if its thickness\footnote{The thickness of a graph $G$ is the minimum number of planar subgraphs whose union equals to $G$.} is 1 or 2.
The {\em complete graph} with $n$ vertices, denoted by $K_n$, is a graph that has an edge between every pair of its vertices. Let $G$ be a subgraph of $K_n$ that has $n$ vertices. The {\em complement} of $G$, denoted by $\overline{G}$, is the graph obtained by removing all edges of $G$ from $K_n$.

As early as 1960 it was known that $K_8$ is biplanar and $K_{11}$ is not biplanar. There exist several biplanar embeddings of $K_8$; see e.g. \cite{Beineke1997} for a self-complementary drawing. The non-biplanarity of $K_{11}$ is easily seen, since it has $55$ edges while a planar graph with eleven vertices cannot have more than $27$ edges, by Euler's formula. Finding the smallest integer $n$, for which $K_n$ is non-biplanar, was a challenging question for some time \cite{Harary1961}.
The following fundamental theorem due to Battle, Harary, and Kodama~(\cite{Battle1962}, 1962) and independently proved by Tutte~(\cite{Tutte1963}, 1963) answers this question and implies that $K_9$ is non-biplanar.

\begin{theorem}
	\label{Harary-Tutte}
	Every planar graph with at least nine vertices has a nonplanar complement.
\end{theorem}

Both proofs of Theorem~\ref{Harary-Tutte} involve a thorough case analysis. Battle, Harary, and Kodama gave an outline of a proof through six propositions. Some of these propositions require detailed case analysis, which is not given in the original paper. For example, the authors write: ``There are several cases to discuss in order to establish Propositions 4 and 5. In each case, we can prove that $\overline{G}$ contains a subgraph homeomorphic to $K_{3,3}$ or $K_{5}$.'' A detailed proof of these propositions is appeared in the master's thesis of Hearon \cite{Hearon2016}.
Tutte's proof is a 13-page paper, and enumerates all simple triangulations (with no separating triangles) with up to 9 vertices and verifies that the complement of each triangulation is nonplanar (the connection to triangulations will become clear shortly).
It seems that Harary was not quite satisfied with any of these proofs as he noted in his Graph Theory book \cite{Harary1969} that ``this result was proved by exhaustion; no elegant or even reasonable proof is known.''
We are still unaware of any short proof of this result. (See \cite{Kuila2014} for a recent attempt towards a new proof.)

The non-biplanarity of $K_9$ has the same flavor as the well-known theorem of Kuratowski on non-planar graphs (stated in Theorem~\ref{Kuratowski-thr}).   The {\em biplanar crossing number} of a graph is the minimum number of crossings over all drawings of the graph in two planes \cite{Czabarka2006}. It is known that $K_9$ can be drawn in two planes with one crossing (see e.g. \cite{Durocher2016}). This and  Theorem~\ref{Harary-Tutte} imply that the biplanar crossing number of $K_9$ is 1. Determining biplanar crossing numbers of $K_n$ for small values of $n$ is important as they lead to better bounds for biplanar crossing numbers of $K_n$ for large values of $n$; see e.g. \cite{Czabarka2006,Czabarka2008,Owen1971}, and  \cite{Durocher2016,Shavali2019} for more recent progress.

\section{Our proof}

\begin{figure}[htb]
	\vspace{-8pt}
	\centering
	\setlength{\tabcolsep}{0in}
	$\begin{tabular}{cc}
		\multicolumn{1}{m{.35\columnwidth}}{\centering\includegraphics[width=.19\columnwidth]{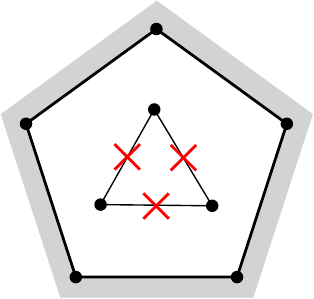}}
		&\multicolumn{1}{m{.6\columnwidth}}{\centering\includegraphics[width=.50\columnwidth]{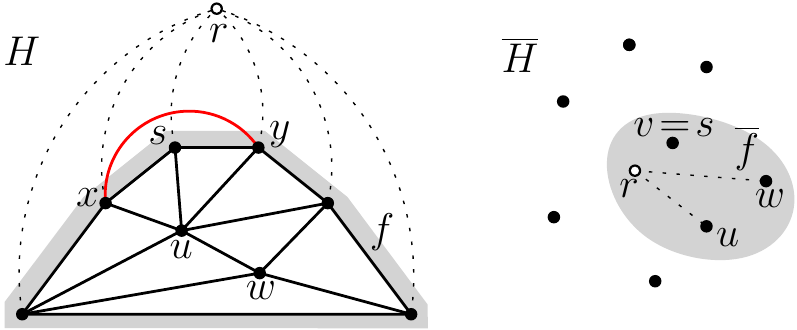}}\\
		(a)&(b)
	\end{tabular}$
	\vspace{-3pt}
	\caption{Illustration of (a) the statement of Theorem~\ref{K8-thr} (b) the proof of Theorem~\ref{Harary-Tutte}.}
	\label{K9-fig}
	\vspace{-12pt}
\end{figure}

In this section we present a short proof of Theorem~\ref{Harary-Tutte}. Our proof is complete, self-contained, and only uses Kuratowski's theorem for non-planar graphs.
Towards our proof we show (in Theorem~\ref{K8-thr}) that a particularly restricted drawing
of $K_8$ cannot be biplanar (see Figure~\ref{K9-fig}(a) for an illustration).
\begin{theorem}
	\label{K8-thr}
	Let $H$ be an embedded planar graph with eight vertices such that the boundary of its outer face is a $5$-cycle and there are no edges between the three vertices that are not on the outer face. Then the complement of $H$ is nonplanar.  
\end{theorem}
\vspace{5pt}

\noindent{\em Proof of Theorem~\ref{Harary-Tutte}.}
Consider a planar graph $G$ with nine vertices. For the sake of contradiction assume that its complement $\overline{G}$ is also planar. Fix a planar embedding of $G$ and a planar embedding of $\overline{G}$. For convenience we use $G$ and $\overline{G}$ for referring to planar graphs and to their planar embeddings. If there are two vertices in $G$ that lie on the same face and are not connected by an edge, then we transfer the corresponding edge from $\overline{G}$ to $G$ and connect the two vertices by a curve in that face.
After this operation both $G$ and $\overline{G}$ remain planar. Repeating this process converts $G$ to an edge-maximal planar graph. In particular $G$ becomes a triangulation in which the boundary of every face (including the outer face) is a triangle (i.e. a 3-cycle).

\vspace{10pt}
\noindent\emph{Claim 1. At least one vertex on the outer face of $G$ has degree larger than four.} To prove this claim we use contradiction. Assume that all three vertices on the outer face of $G$ are of degree at most 4. The removal of these three vertices from $G$ results in a $6$-vertex graph $G'$. The region, that is between the boundaries \begin{wrapfigure}{r}{1.2in} 
	\centering
	\vspace{-15pt}
	\includegraphics[width=1in]{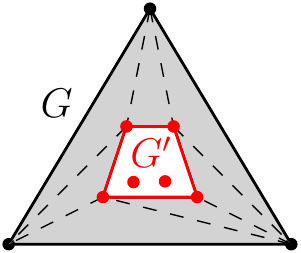} 
	\vspace{-10pt} 
\end{wrapfigure}of the outer face of $G$ and the outer face of $G'$ is a polygon with a hole, that is triangulated by at most six edges of $G$ (because every vertex on the outer face of $G$ has at most two edges in the interior of this polygon). The boundary of the outer face of $G'$, i.e. the hole, has three vertices because otherwise (if it has at least four vertices) the polygon would require at least seven edges to be triangulated, as in the figure to the right; this can be verified by a simple counting argument using Euler's formula for planar graphs, see also \cite[Proof of Lemma 5.2]{Rourke1987}. Thus  the outer face of $G'$ is a $3$-cycle. In this case the other three vertices of $G'$ which are in the interior of this 3-cycle together with the three removed vertices from $G$ form a $K_{3,3}$ in $\overline{G}$, which contradicts its planarity. This proves Claim 1.

\vspace{7pt}
In view of Claim 1 we assume that at least one vertex, say $r$, on the outer face of $G$ has degree $k\geq 5$. Remove $r$ from $G$ and $\overline{G}$ and denote the resulting graphs by $H$ and $\overline{H}$, respectively. Notice that $(H,\overline{H})$ is a biplanar embedding of $K_8$. Let $f$ and $\overline{f}$ be the faces of $H$ and $\overline{H}$, respectively, that contain the removed vertex $r$, as in Figure~\ref{K9-fig}(b). Notice that $f$ is the outer face of $H$. Since $(G,\overline{G})$ was a biplanar embedding of $K_9$, in which $r$ was connected to all other 8 vertices, we have the following observation.

\vspace{8pt}
\noindent{\em Observation 1. Every vertex of the resulting graph $K_8$ lies on $f$ or on $\overline{f}$.}

Since $\overline{G}$ was a simple graph (no multiedges and no loops), the face $\overline{f}$ has at least three vertices; these vertices are not necessarily connected in $\overline{H}$. Since $G$ was a triangulation, the boundary of the outer face $f$ of $H$ is a $k$-cycle. 
If $k>5$ then let $s$ be a vertex of $f$ that also lies on $\overline{f}$; such a vertex exists because $\overline{f}$ has at least three vertices and we have eight vertices in total. Let $x$ and $y$ be the neighbors of $s$ on $f$. If $xy$ is an edge of $H$ then draw it as a curve in $f$. If $xy$ is not an edge of $H$ then transfer it from $\overline{H}$ to $H$ and draw it in $f$, as in Figure~\ref{K9-fig}(b). Now, the new outer face $f$ of $H$ has $k-1$ vertices. Repeat the above process until the outer face of $H$ has exactly five vertices. 

At this point $f$ has five vertices. 
Let $u,v,w$ be the vertices of $K_8$ that are
not on $f$. These three vertices lie on $\overline{f}$, because of Observation 1 and our choices of $s$ (for the case $k>5$). If any of the edges $uv$, $uw$, and $vw$ are not in $\overline{H}$ then transfer them from $H$ to $\overline{H}$ and draw in $\overline{f}$ without crossing other edges. We obtain a planar graph $H$ that satisfies the constraints of Theorem~\ref{K8-thr} and so that its complement $\overline{H}$ is planar. This contradicts Theorem~\ref{K8-thr}.~\hfill$\qed$
\vspace{8pt}

To prove Theorem~\ref{K8-thr} we use the theorem of Kuratowski  \cite{Dirac1954,Kuratowski1930} that ``a finite graph is non-planar if and only if it contains a subgraph that is homeomorphic to $K_5$ or $K_{3,3}$.'' The following is an alternative statement for Kuratowski's theorem, which is given in \cite{Tutte1963}.

\begin{theorem}
	\label{Kuratowski-thr}
	A graph $G$ is nonplanar if one of the following conditions hold: 
	$(i)$ $G$ has six disjoint connected subgraphs $A_1, A_2, A_3, B_1, B_2, B_3$ such that for each $A_i$ and $B_j$ there is an edge with one end in $A_i$ and the other in $B_j$. $(ii)$ $G$ has five disjoint connected subgraphs $A_1, A_2, A_3, A_4, A_5$ such that for each $A_i$ and $A_j$, with $i\neq j$, there is an edge with one end in $A_i$ and the other in $A_j$.
\end{theorem}

\vspace{5pt}
{\noindent \em Proof of Theorem~\ref{K8-thr}.} 
Let the $5$-cycle $C=(a_1,a_2,a_3,a_4,a_5)$ be the boundary of the outer face of $H$, and let $u$, $v$, and $w$ be the three vertices that are not on the outer face, i.e., lie on internal faces of $H$. By the statement of the theorem $uv$, $uw$, and $vw$ are edges of the complement graph $\overline{H}$.
Except for the three pairs $(u,v)$, $(u,w)$, $(v,w)$, if a pair of vertices lie on the same internal face of $H$ and are not connected by an edge, then we transfer the corresponding edge from $\overline{H}$ to $H$ and connect the two vertices by a curve in the face. 
After this operation $H$ remains planar. Repeating this process makes $H$ edge-maximal (in the above sense). 

\begin{figure}[htb]
	\centering
	\setlength{\tabcolsep}{0in}
	$\begin{tabular}{ccc}
		\multicolumn{1}{m{.33\columnwidth}}{\centering\includegraphics[width=.25\columnwidth]{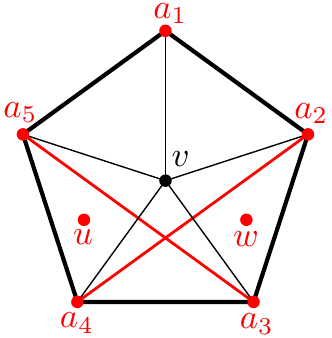}}
		&\multicolumn{1}{m{.33\columnwidth}}{\centering\includegraphics[width=.25\columnwidth]{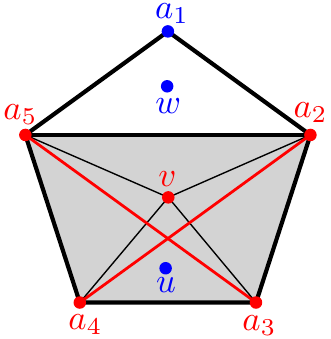}}
		&\multicolumn{1}{m{.33\columnwidth}}{\centering\includegraphics[width=.25\columnwidth]{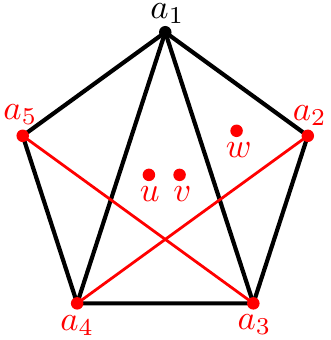}}
		\\
		(a)&(b)&(c)
	\end{tabular}$
	\caption{Black edges belong to $H$, bold black edges belong to $H'$, and red edges belong to $\overline{H}$.}
	\label{K8-fig}
\end{figure}

Let $H'$ be the embedded planar subgraph of $H$ that is induced by the five vertices of $C$. The graph $H'$ consists of the cycle $C$ together with zero, one, or two chords as in Figure~\ref{K8-fig}.

\vspace{8pt}

\begin{wrapfigure}{r}{1.15in} 
	\centering
	\vspace{-18pt} 
	\includegraphics[width=1.05in]{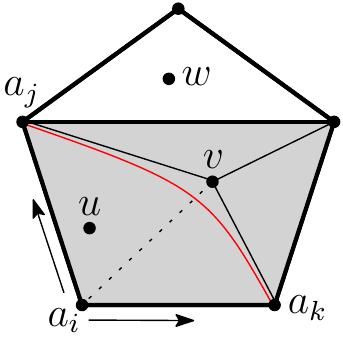} 
	\vspace{-8pt} 
\end{wrapfigure}
\noindent\emph{Claim 2. If an internal face $f$ of $H'$ contains $u$, $v$, or $w$ then one of them is connected to all boundary vertices of $f$ in $H$.} The shaded region in the figure to the right represents $f$. To verify the claim, first observe that (by edge-maximality of $H$) one of the vertices in $f$, say $v$, is connected to at least three boundary vertices of $f$, i.e., $v$'s degree in $H$ is at least three. We argue that $v$ should be connected to all boundary vertices of $f$. For a contradiction assume that $v$ is not connected to some vertex $a_i$ on $f$. Let $a_j$ and $a_k$ be the neighbors of $v$ on $f$ that are visited first while walking on boundary of $f$ in clockwise and counterclockwise directions starting from $a_i$. Since $v$ is not connected to other vertices in the interior of $f$, we could have moved the edge $a_ja_k$ from $\overline{H}$ to $H$ and draw it in $f$. This means that $H$ is not edge-maximal, which is a contradiction.
\vspace{8pt}

Now we consider three cases depending on the number of chords of $H'$. In each case we show that $\overline{H}$ is nonplanar. 
\begin{itemize}
	\item $H'$ has no chords. Let $v$ be the vertex of $H$ that (by Claim 2) is connected to each $a_i$; see Figure~\ref{K8-fig}(a).
	By planarity of $H$, each of $u$ and $w$ can only be adjacent to two consecutive vertices of $C$. Hence there
	exists a vertex of $C$ (say $a_1$) that is adjacent
	to neither $u$ nor $w$. 
	In this setting, regardless of the locations of $u$ and $w$, the five connected subgraphs $u$, $w$, $a_1$, $\{a_2,a_4\}$ and $\{a_3,a_5\}$ from $\overline{H}$ satisfy condition (ii) of Theorem~\ref{Kuratowski-thr}. Thus $\overline{H}$ is nonplanar. 
	\item $H'$ has one chord. After a suitable relabeling assume that this chord is $(a_2,a_5)$. Let $f$ denote the face of $H'$ whose boundary is the 4-cycle $(a_2,a_3,a_4,\allowbreak a_5)$; this face is shaded in Figure~\ref{K8-fig}(b). This face contains some vertices of $\{u,v,w\}$ because otherwise $H'$ should have a chord in $f$ (by maximality of $H$) which contradicts our assumption that $H'$ has one chord. Let $v$ be the vertex in $f$ that (by Claim 2) is connected to all its boundary vertices. By planarity of $H$, each of $u$ and $w$ can only be adjacent to two consecutive vertices of $f$. Therefore, the six connected subgraphs $u$, $w$, $a_1$, $v$, $\{a_2,a_4\}$, and $\{a_3,a_5\}$ from $\overline{H}$ (partitioned into $\{u, w, a_1\}$ and $\{v, \{a_2,a_4\},\{a_3,a_5\}\}$) satisfy condition (i) of Theorem~\ref{Kuratowski-thr}. Thus $\overline{H}$ is nonplanar.
	\item $H'$ has two chords. Let $a_1$ be the vertex that is incident to the two chords as in Figure~\ref{K8-fig}(c). By planarity of $H$, each of $u$, $v$, and $w$ can only be adjacent to one vertex in $\{a_2,a_4\}$ and to one vertex in $\{a_3,a_5\}$. Thus, the five connected subgraphs $u$, $v$, $w$, $\{a_2,a_4\}$, and $\{a_3,a_5\}$ from $\overline{H}$ satisfy condition (ii) of Theorem~\ref{Kuratowski-thr}, and hence $\overline{H}$ is nonplanar.~\hfill$\qed$
\end{itemize}

\vspace{-15pt}
\vspace{-25pt}
%
%
%
%

\section{Conclusions}
For any integer $k\ge 1$ let $\nu(k)$ be the smallest integer for which the (edges of the) complete graph with $\nu(k)$ vertices cannot be drawn in $k$ planes without creating a crossing. As the maximum number of (noncrossing) edges that can be drawn in a plane is $3\nu(k)-6$ and the number of edges of the complete graph is $\nu(k)\choose 2$, a counting argument implies that 
$$\nu(k)\le \left\lfloor\frac{6k+1+\sqrt{36k^2-36k+1}}{2}\right\rfloor+1.$$   
This bound implies that $\nu(1)\le 5$ and $\nu(2)\le 11$, however for $k\in\{1,2\}$ we already know that $\nu(1)=5$ and $\nu(2)=9$. It would be interesting to find exact value of $\nu(k)$ for larger values of $k$.
\bibliographystyle{splncs04}
\bibliography{Non-biplanar-K9.bib}
\end{document}